\documentstyle{article}

\def\ni{\noindent }
\def\eq #1{(\ref{#1})}       
\def\l{\left}                   
\def\r{\right}                  

\newcommand{\be}[1]{\begin{equation}\label{#1}}
\def\ee{\end{equation}}

\newcommand{\ba}[1]{\begin{array}{#1}}
\def\ea{\end{array}}

\def\fr #1#2{\frac{#1}{#2}}

\def\se #1{sec.\,\ref{#1}}

\def\y1{\mbox{$y'$}}
\def\yt{\mbox{$y''$}}
\def\yyy{\mbox{$y'''$}}

\def\hyper3{{\bf hyper3}}   

\def\3F2{\mbox{$_3${F}$_2$}}
\def\2F2{\mbox{$_2${F}$_2$}}
\def\1F2{\mbox{$_1${F}$_2$}}
\def\0F2{\mbox{$_0${F}$_2$}}
\def\pFq{\mbox{$_p${F}$_q\;$}}

\def\PFQ{\mbox{$_p${F}$_q\;$}}
\begin{document}
\title{Hypergeometric solutions\\ for third order linear ODEs}

\author{E.S. Cheb-Terrab$^{a}$ and A.D. Roche $^{a,b}$}
\date{}
\maketitle
\thispagestyle{empty}


\medskip
\centerline {\it $^a$Maplesoft, Waterloo Maple Inc.}

\medskip
\centerline {\it $^b$Department of Mathematics}
\centerline {\it Simon Fraser University, Vancouver, Canada.}


\maketitle

\begin{abstract}


In this paper we present a decision procedure for computing \pFq hypergeometric solutions for third order linear ODEs, that is, solutions for
the classes of hypergeometric equations constructed from the \3F2, \2F2, \1F2
and \0F2 standard equations using transformations of the form $x \rightarrow
F(x),\, y \rightarrow P(x) y$, where $F(x)$ is rational in $x$ and $P(x)$ is
arbitrary. A computer algebra implementation of this work is present in
Maple 12.

\end{abstract}

\section*{Introduction}

Given a third order linear ODE
\be{LinearODE}
y''' + c_2\,y'' + c_1\,y' + c_0\,y = 0
\ee

\ni where $y \equiv y(x)$ is the dependent variable and the $c_j \equiv
c_j(x)$ are any functions of $x$ such that the quantities\footnote{$I_1$ and $I_0$ are invariant under transformations of the dependent variable of
the form $y(x) \rightarrow P(x)\,y(x)$, $P$ arbitrary.}

\be{invariants}
I_1 = {c_{{2}}}^{'} + \fr{{c_{{2}}}^{2}}{3} - c_{{1}}\,
\ \ \ I_0 = \fr{{c_{{2}}}^{{\it ''}}}{3} - {\frac{2\,{c_{{2}}}^{3}}{27}} + \fr{c_{{1}}c_{{2}}}{3} - c_{{0}}
\ee

\ni are rational functions of $x$, the problem under consideration is that of
systematically computing solutions for \eq{LinearODE} even when no Liouvillian solutions
exist\footnote{Expressions that can be expressed in terms of exponentials,
integrals and algebraic functions, are called Liouvillian. The typical example
is $\exp(\int R(x), dx)$ where $R(x)$ is rational or an algebraic function
representing the roots of a polynomial.}. Recalling, Liouvillian
solutions can be computed systematically
\cite{LiouvillianSolutions} and implementations of the related algorithm exist
in various computer algebra systems. The linear ODEs involved in mathematical
physics formulations, however, frequently admit only non-Liouvillian special
function solutions, and for this case the existing algorithms cover a rather
restricted portion of the problem.

The special functions associated with linear ODEs frequently happen to be
particular cases of some generalized hypergeometric \PFQ functions \cite{Seaborn}. One
natural approach is thus to directly search for \pFq solutions instead of
special function solutions of one or another kind, and this is the approach
discussed here. Related computer algebra routines were implemented in 2007 and
are now at the root of the Maple (release 12) \cite{maple12} ability for
solving non-trivial 3rd order linear ODE problems.

The approach used consists of resolving an equivalence problem between a given
equation of the form \eq{LinearODE} and the four standard \PFQ differential
equations associated to third order linear ODEs, that is, the \3F2, \2F2, \1F2
and \0F2 equations \cite{abramowitz}, respectively:

\be{seeds}
\begin{array}{l}
\displaystyle
\yyy - {\frac { \l(  \delta+\eta+1 - \l( \alpha+\beta+\gamma+3 \r) x \r) }{x \l( x-1 \r) }}\,\yt
	- {\frac { \l(  \eta\,\delta - \l(  \l( \beta+\gamma+1 \r) \alpha + \l( \beta+1 \r) \l( \gamma+1 \r)  \r) x \r)}{{x}^{2} \l( x-1 \r) }}\, \y1
	+ {\frac {\alpha\,\beta\,\gamma}{{x}^{2} \l( x-1 \r) }}\,y = 0
\\*[0.12in]
\displaystyle
\yyy - {\frac { \l( x-\gamma-\delta-1 \r) }{x}}\,\yt
	- {\frac { \l(  \l( \alpha+\beta+1 \r) x-\gamma\,\delta \r) }{{x}^{2}}}\,\y1
	- {\frac {\alpha\,\beta}{{x}^{2}}}\,y = 0
\\*[0.12in]
\displaystyle
\yyy + {\frac { \l( \beta+\gamma+1 \r) }{x}}\,\yt
	- {\frac { \l( x-\beta\,\gamma \r) }{{x}^{2}}}\,\y1
	- {\frac {\alpha}{{x}^{2}}}\,y = 0
\\*[0.12in]
\displaystyle
\yyy + {\frac { \l( \alpha+\beta+1 \r) }{x}}\,\yt
	+ {\frac {\alpha\,\beta}{{x}^{2}}}\,\y1
	- {\frac {1}{{x}^{2}}}\,y = 0
\end{array}
\ee
\smallskip

\ni where $\{\alpha, \beta, \gamma, \delta, \eta\}$ represent arbitrary
expressions constant with respect to $x$. The equivalence classes are
constructed by applying to these equations the general transformation\footnote{The
problem of equivalence under transformations $\{x \rightarrow F(x),\ \ y
\rightarrow P(x)\, y + Q(x)\}$ for linear ODEs can always be mapped into one
with $Q(x) = 0$, see \cite{ince}.}

\be{tr}
x \rightarrow F(x),\, y \rightarrow P(x)\, y
\ee

\ni where $P(x)$ is arbitrary, with the only restriction that $F(x)$ is
rational in $x$, resulting in rather general ODE families. When the equation being solved belongs to this class, apart from
providing the values of $F(x)$ and $P(x)$ that resolve the problem, the
algorithm systematically returns the values of the (five, four, three or
two) \pFq parameters entering each of the three independent solutions.

It is important to note that the idea of seeking hypergeometric function
solutions for linear ODEs or using an equivalence approach for that purpose is
not new, although in most cases the approaches presented only handle second
order linear equations \cite{kamran, willis, bronstein, hyper3}. An
exception to that situation are the algorithms \cite{petkovsek, george_MeijerG} for
computing \pFq solutions for third and higher order linear ODEs, and a similar
one implemented in Mathematica \cite{mathematica}. It is our understanding, however, that the transformations defining the classes of
equivalence that those algorithms can handle are restricted to
$x \rightarrow a x^b,\, y \rightarrow P(x) y$, with $a$ and $b$ constants, not
having the generality of \eq{tr} with rational $F(x)$ presented here.

Apart from concretely expanding the ability to solve third order linear
ODEs, the decision procedure being presented generalizes previous work in that:

\begin{enumerate}

\item The ideas presented in \cite{hyper3}, useful for
decomposing two sets of invariants into each other,
were extended for third order equations and elaborated further.

\item The classification ideas presented in \cite{hyper3} for second order
linear equations were extended for third order.

\item When the \pFq parameters are such that less than three independent \pFq
solutions exist, instead of introducing integrals \cite{george_MeijerG}, MeijerG
functions are used to express the missing independent solutions.

\end{enumerate}

The combination of items 1 and 2 resulted in the new ability to solve the \pFq ODE
classes generated by transformations as general as \eq{tr} with $F(x)$ rational. Item 3
is not new\footnote{Mathematica 6 also uses MeijerG functions as described in
item 3.}, 
though we are not aware of literature presenting the related problem
and solution. Altogether, these ideas and its related algorithm permit
the systematic computation of three
independent solutions for a large set of third order linear equations that
we didn't know how to solve before.


\section{Computing hypergeometric solutions}
\label{hyper3}

To compute \PFQ solutions to \eq{LinearODE} the idea is to formulate an
equivalence approach to the underlying hypergeometric
differential equations, that is, to determine whether a given linear ODE can
be obtained from one of the \PFQ ODEs \eq{seeds} by means of a transformation
of a certain type. If so, the solution to the given ODE is obtained
by applying the same transformation to the solution of the corresponding
\PFQ equation.


The approach also requires determining the values of the
hypergeometric parameters $\{\alpha, \beta, \gamma, \delta, \eta\}$ for which
the equivalence exists, and it is clear that the bottleneck in this approach
is the generality of the class of transformations to be considered. For
instance, one can verify that for linear transformations of the form \eq{tr}
with arbitrary $F(x)$, in the case of second order linear ODEs, the problem
is too general in that the determination of $F(x)$ requires solving the given
ODE itself \cite{kathi}, making the approach of no practical use. This has to
do with the fact that in the second order case, any linear ODE can be
obtained from any other one through a transformation of the form \eq{tr}. The
situation for third order equations is different: the transformation \eq{tr}
is not enough to map any equation into any other one\footnote{Therefore there
exist enough absolute invariants to formulate the equivalence problem under
\eq{tr} - see \se{MinimalDegrees}.} \cite{MSW_4ODE_challenges}, so that its
determination when the equivalence exists is in principle possible. By
restricting the form of $F(x)$ entering \eq{tr} to be rational in $x$ the
problem becomes tractable by using a two step strategy:

\begin{enumerate}

\item Compute a rational transformation $R(x)$ mapping the normal form of the
given equation\footnote{The coefficients of $\y1$ and $y$ in the
normalized equation are the invariants $I_1$ and $I_0$ defined in
\eq{invariants}, assumed to be rational.} into one having {\it invariants with minimal degrees}
(defined in \se{MinimalDegrees}).

\item Resolve an equivalence problem between this equation with minimal
degrees and the standard $_p\mbox{F}_q$ equations \eq{seeds} under 
transformations of the form discussed in \cite{hyper3}, that is

\be{tr_hyper3}
x \rightarrow \fr{(a\,x^k+b)}{(c\,x^k+d)},\ \ \ \ y \rightarrow P(x)\,y
\ee

\ni with $P(x)$ arbitrary and $\{a, b, c, d, k\}$ constants with respect to
$x$. In doing so, determine also the parameters $\{\alpha, \beta, \gamma,
\delta, \eta\}$ of the \pFq or MeijerG functions entering the three
independent solutions.

\end{enumerate}

\ni The key observation in this ``two steps" approach is that a
transformation of the form \eq{tr} with rational $F(x)$ mapping into the \pFq equations \eq{seeds}, when it exists, it can always be
expressed as the composition of two transformations, each one related to each
of the two steps above (see \se{MinimalDegrees}), because \eq{seeds} have invariants with minimal degrees. The advantage of splitting the problem in this way is
that the determination of $R(x)$ in step one, and of the
(up to five) \pFq parameters in step two, as well as of the values of $\{a,
b, c, d, k\}$ entering \eq{tr_hyper3}, is systematic (see \se{equiv_section} and \se{MinimalDegrees}), even when the
problem is nonlinear in many variables.

\section{Equivalence under $x \rightarrow {(a\,x^k+b)}/{(c\,x^k+d)},\  y \rightarrow P(x)\,y$}
\label{equiv_section}

This type of equivalence is discussed in \cite{hyper3} and
generalized here for third order ODEs. Recalling the main points, these
transformations, which do not form a group in the strict sense, can be
obtained by sequentially composing three different transformations, each of
which does constitute a group. The sequence starts with linear fractional -
also called M\"obius - transformations

\be{mobius}
x \rightarrow \fr{a\, x + b}{c\, x + d},
\ee

\ni is followed by power transformations

\be{power}
x \rightarrow x^k,
\ee

\ni and ends with linear homogeneous transformations of the dependent variable

\be{linear_in_y}
y \rightarrow P\,y.
\ee

\subsection{Equivalence under transformations of the dependent variable $y \rightarrow P(x) y$}

Transformations of the form \eq{linear_in_y} can easily be factored out of
the problem: if two equations of the form \eq{LinearODE} can be obtained from
each other by means of \eq{linear_in_y}, the transformation relating them is
computable directly from these coefficients. For that purpose first rewrite
both equations in normal form using

\be{ToNormalForm}
\displaystyle
y \rightarrow  y\,{e^{-\int \!c_2(x)/3\,{dx}}}
\ee

\ni and the transformation relating the two hypothetical ODEs - say with
coefficients $c_j$ and $\tilde{c}_k$, when it exists, is given by $y
\rightarrow y\,{e^{\int \!(c_2(x)-\tilde{c}_2(x))/3\,{dx}}}$.

\subsection{Equivalence under M\"obius transformations, singularities and classification}
\label{equivalence_under_mobius}

M\"obius transformations preserve the structure of the singularities of \eq{LinearODE}. For example, all of the
\0F2, \1F2 and \2F2 hypergeometric equations in \eq{seeds} have one regular
singularity at the origin and one irregular singularity at infinity, and
after transforming them using the M\"obius transformations \eq{mobius}, they
continue having one regular singularity and one irregular singularity, now
respectively located at\footnote{When either $a$ or $c$ are equal to zero,
the corresponding singularity is located at $\infty$} $-b/a$ and $-d/c$.

In the case of the \3F2 differential equation (the first listed in
\eq{seeds}), under \eq{mobius} the three regular singularities move from
$\{0,1,\infty\}$ to $\{-b/a,\, -d/c,\, (d-b)/(a-c)\}$. So from the
singularities of an ODE, not only one can tell with respect to which of the
four differential equations \eq{seeds} could the equivalence under
\eq{mobius} be resolved, but also one can extract the values of the
parameters $\{a, b, c, d\}$ entering the transformation \eq{mobius}.

More generally, through M\"obius transformations one can formulate a 
classification of singularities of the linear ODEs ``equivalent" to the third
order \PFQ equations \eq{seeds} as done in \cite{hyper3}
for second order \PFQ equations. So, for each \pFq family obtained from \eq{seeds} using \eq{mobius}, a
classification table can be constructed based only on:
\begin{itemize}

\item the degrees of the numerators and denominators of the invariants \eq{invariants};

\item the presence of roots with multiplicity in the denominators;

\item the possible cancellation of factors between the numerator and
denominator of each invariant.

\end{itemize}

\ni
With this classification in hands, from the
knowledge of the degrees with respect to $x$ of the numerator and denominator
of the invariants \eq{invariants} of a given third order linear ODE, one can determine systematically
whether or not the equation could be obtained from the \3F2, \2F2, \1F2 or \0F2
equations \eq{seeds} using \eq{mobius}.

\subsection{Transformations $x \rightarrow F(x)$ and equivalence under $x \rightarrow x^k$}

Changing $x \rightarrow F(x)$ in
\eq{LinearODE}, the new invariants $\tilde{I}_j$ can be
expressed in terms of the invariants \eq{invariants} of \eq{LinearODE} by

\be{Inv_1}
\begin{array}{l}
{\tilde I}_{{1}} ( x ) ={{\it F'}}^{2}I_{{1}} ( F ) - 2\,S ( F ) 
\\*[0.15in]
\displaystyle
{\tilde I}_{{0}} ( x ) = {\it F'}{\it F''}I_{{1}} ( F ) +{{\it F'}}^{3}I_{{0}} (F) - S(F)'
\end{array}
\ee

\ni where $S(F)$ is the Schwarzian \cite{weisstein}
\be{Schwarzian}
S(F) = \fr{F'''}{F'} - \fr{3}{2} \l( \fr{F''}{F'} \r)^2.
\ee

\ni The form of $S(F)$ is particularly simple when $F(x)$ is a M\"obius
transformation, in which case $S(F) = 0$.

Regarding power transformations $F(x) = x^k$, 
unlike M\"obius transformations, they {\em do not preserve} the structure of
singularities; the Schwarzian \eq{Schwarzian} is:

\be{S_for_powers}
S(x^k) = {\frac {1 - {k}^{2}}{2\, {x}^{2}}}.
\ee

\ni From \eq{Inv_1} and \eq{S_for_powers}, for
instance the transformation rule for $I_1(x)$ becomes

\begin{equation}
{x}^{2}{\tilde I}_{{1}} ( x ) + 1 = k^2 \l(  ({x}^{k}) ^{2}I_{{1}}( {x}^{k} )+ 1 \r).
\ee

\ni Generalizing to third order the presentation of shifted invariants in \cite{hyper3},
we define here

\be{J_definition}
\begin{array}{l}
J_1(x) = {x}^{2}{I}_{{1}} ( x ) +1,
\\*[0.15in]
J_2(x) = {x}^{3}{I}_{{0}} ( x ) + {x}^{2}{I}_{{1}} ( x ).
\end{array}
\ee

\ni From \eq{Inv_1} rewritten in terms of these $J_n(x)$, their transformation rule
under $x \rightarrow x^k$ is given by

\be{J_transformation_rule}
\begin{array}{l}
{\tilde J}_1(x) = k^2 J_1(x^k),
\\*[0.15in]
{\tilde J}_2(x) = k^3 J_2(x^k).
\end{array}
\ee

The equivalence of two linear ODEs A and B under $x \rightarrow x^k$ can then
be formulated as follows: Given the shifted invariants ${\tilde J}_{n,A}(x)$
and ${\tilde J}_{n,B}(x)$, computed using their definition \eq{J_definition} in terms of
${\tilde I}_n(x)$ defined in \eq{invariants}, compute $k_A$ and $k_B$ entering
\eq{J_transformation_rule} such that the degrees of
${J}_{n, A}(x)$ and ${J}_{n, B}(x)$ are minimal. From the knowledge of $x \rightarrow x^{k_A}$ and $x \rightarrow
x^{k_B}$, respectively leading to $J_{n, A}$ and $J_{n, B}$ with minimized
degrees, equations A and B are related through power transformations only when
$J_{n, A} = J_{n, B}$ and, if so, the mapping relating A and B is $x
\rightarrow x^{k_A-k_B}$. Finally, the computation of $k$ simultaneously
minimizing the degrees of the two $J_n(x)$ in \eq{J_transformation_rule} is
performed as explained in section 3 of \cite{hyper3}.

\section{Mapping into equations having invariants with minimal degrees}
\label{MinimalDegrees}

The decision procedure presented in the previous section serves for systematically solving well defined families
of \pFq 3rd order equations for which no solving algorithm was available
before to the best of our knowledge. However, the restriction in the form of $F(x)$ entering \eq{tr} to the
composition of M\"obius with power transformations is unsatisfactory: for
linear equations of order higher than two, \eq{tr} does not map any linear
equation into any other one of the same order and so the problem is already
restricted\footnote{The equivalence problem for linear equations of order $n$
involves a system of $n-1$ equations and invariants $I_j(x)$, that includes the equivalence
function $F(x)$. When $n > 2$, eliminating $F(x)$ from the problem results in
an interrelation between the $I_j$ so that the equivalence is only possible
when these relationships between the $I_j$ hold \cite{MSW_4ODE_challenges}.}.

As shown in what follows, one possible extension of the algorithm is thus to consider the general
transformations \eq{tr} restricting $F(x)$ to be a rational
function of $x$. For that purpose, instead of working with invariants $I_j$
under $y \rightarrow P(x)\,y$ we introduce absolute invariants $L_i$
under $\{x \rightarrow F(x),\  y \rightarrow P(x)\,y\}$:

\be{Ls}
L_1 = \fr{(6rr''+9 I_1 r^2-7{r'}^2)^3}{r^8},\ \ \ 
L_2 = \fr{(27 {I_1}' r^3-18 I_1 r^2 r'+56{r'}^3-72 r'' r' r+18 r''' r^2)}{r^4};
\ee

\ni where $r={I_1}' - 2 I_0$ is a relative invariant of weight 3 \cite{wilczynski}.
Under \eq{tr}, $L_i$ transforms as $L_i(x) \to L_i(F(x))$ and \eq{Ls}
can be inverted using as intermediate variables the relative
invariants $s=(L_2 L_1)/{L_1}'$, and $t=L_1/s^3$:

\be{Is}
I_1 = \fr{st^3-6t''t+7{t'}^2}{9t^2},\ \ \ 
I_0 = \fr{(s'-9)t^4+t'st^3-6t'''t^2+20t''t't-14{t'}^3}{18 t^3}.
\ee

Thus, any canonical form for the $L_i$ that can be achieved using \eq{tr}
automatically implies on a canonical form for the $I_i$ and so for the ODE
\eq{LinearODE}. The canonical form we propose here is one where the $L_i$ have
{\it minimal degrees}, that is, where the maximum of the degrees of the
numerator and denominator in each of the $L_i(x)$ is the minimal one that can
be obtained using a rational transformation $x \to F(x)$. This canonical form
is not unique in that it is still possible to perform a M\"obius
transformation \eq{mobius}, that changes the $L_i$ but not their degrees.

The equivalence of two linear ODEs A and B under \eq{tr} with rational $F(x)$
can then be formulated by rewriting both equations in this canonical form,
where the invariants $L_i$ of each equation have minimal degrees, followed by
determining whether these canonical forms are related through a M\"obius
transformation.

In the framework of this paper, B is one of the hypergeometric equations
\eq{seeds}, all of them already in canonical form in that the corresponding
$L_i$ already have minimal degrees. Hence, the equivalence of A, of the form
\eq{LinearODE}, and any of B of the form \eq{seeds} requires determining only
a canonical form for A (the rational function $F(x)$ minimizing the degrees of
the $L_i$ of A), followed by resolving an equivalence under M\"obius
transformations between this canonical form and any of the equations
\eq{seeds}, done as explained in \se{equivalence_under_mobius}.

The key computation in this formulation of the equivalence problem under
\eq{tr} is thus the computation of a rational $F(x)$ that minimizes the
degrees of the $L_i$ of A. The computation of $F(x)$ can clearly be formulated
as a rational function decomposition problem subject to constraints: {\em
``given two rational functions $L_i(x), i = 1..2$, find rational functions
$\tilde{L}_i(x)$ and $F(x)$ satisfying $L_i = \tilde{L}_i \circ F$ and such
that the rational degree of $F$ is maximized"} (and therefore the degrees of
the canonical invariants $\tilde{L}_i$ are minimized). In turn, this
type of function decomposition associated to ``minimizing the
degrees" of the $L_i$ can be interpreted as the reparametrization, in terms
of polynomials of lower degree, of a rational curve that is improperly
parameterized, as discussed in \cite{sedenberg}, where an algorithm to perform
this reparametrization is presented.

One key feature of the algorithm presented in \cite{sedenberg} is that it reduces
the computation of $F(x)$ to a sequence of {\em univariate} GCD computations,
avoiding the expensive computation of bivariate GCD. However, it is not clear
for us whether the prescriptions in \cite{sedenberg} (at page 71)
for mapping the bivariate GCDs into univariate ones is complete. We also
failed in obtaining a copy of the computer algebra packages {\em FRAC}
\cite{frac} or {\em Cadecom} \cite{cadecom} that contain an implementation of the algorithm presented in
\cite{sedenberg}. Mainly for these reasons, and without the intention of being
original, we describe here a slightly modified version of the algorithm
presented in \cite{sedenberg}.

\subsection{An algorithm for computing $x \rightarrow F(x)$ minimizing the degrees of the $L_i(x)$}

Let $L_i(x) = \tilde{L}_i(F(x)) = N_i(x)/D_i(x), i=1..n$, and
$F(x)=p(x)/q(x)$, where the $\tilde{L}_i$ have minimal degrees, $N_i$ is
relatively prime to $D_i$ and $p$ is relatively prime to $q$.
Construct polynomials 

\be{}
Q_i(x,t) = \mbox{numerator}(L_i(x)-L_i(t)) = N_i(x)D_i(t)-N_i(t)D_i(x),
\ee

\ni and let $P(x,t)$ be the bivariate GCD of these $Q_i(x,t)$. Consequently

\be{P_i}
P(x,t) = \mbox{numerator}(F(x)-F(t)) = \sum_i P_i(x)\, t^i = p(x)\,q(t)-p(t)\,q(x),
\ee

\ni The coefficient $P_i(x)$ of each power of $t$ in $P(x,t)$
is a linear combination of $p(x)$ and $q(x)$, and because the quotient of any
two relatively prime of these linear combinations is fractional linear
in $F(x)$, so is the quotient of any two relatively prime $P_i(x)$. Finally,
because $F(x)$ is defined up to a M\"obius transformation we can take that quotient itself - say, $P_i(x)/P_j(x)$ - as the solution $F(x)$.

The slowest step of this algorithm is the computation of the bivariate GCD
between the $Q_i(x, t)$ that determines the function $P(x,t)$ from which the $P_i(x)$ are computed. It is possible
however to avoid computing that
bivariate GCD, using a small number of univariate GCD computations instead.

For that purpose, notice first that what is relevant in the $P_i(x)$ is that they are
linear combinations of $p(x)$ and $q(x)$.
Now, we can also obtain linear combinations of $p(x)$ and $q(x)$ by directly
substituting numerical values $t_k$ for $t$ into $P(x,t)$,
and from there compute $F(x)$ as the quotient, e.g., of $P(x,t_0)/P(x,t_1)$.
The key observation here
is that these $P(x, t_k)$ can also be obtained by
substituting $t=t_k$ directly into the $Q_i(x,t)$ followed by computing the univariate GCD of
$Q_1(x,t_k)$ and $Q_2(x,t_k)$\footnote{For example, suppose the $x$-solutions of
$P(x,t)=0$ are $x=X_j(t), j=1..m$, i.e., $P(x,t) = -P_m(t) \prod_j x-X_j(t)$.
Then each $x=X_j(t_0)$ is a solution of both $Q_1(x,t_0)=0$ and
$Q_2(x,t_0)=0$. For most values of $t_0$ (all but a finite set in fact) these
$X_j(t_0)$ will be the only such common solutions, and therefore the GCD of
$Q_1(x,t_0)$ and $Q_2(x,t_0)$ is in fact $P(x,t_0)$.}, avoiding in this way the computation of the expensive bivariate GCD leading to $P(x,t)$.


Repeating this process with another $t$-value gives a second, in general
different, such linear combination of $p(x)$ and $q(x)$, with $F$ being the resulting quotient of two of these linear combinations obtained using different values of $t$.
The rest of the algorithm entails avoiding invalid $t$-values at the time of
substituting $t = t_k$ and this is accomplished by considering different $t_k$
until the following conditions are both satisfied:

\begin{enumerate}

\item The two $P(x,t_0)$, $P(x,t_1)$ whose quotient gives the solution $F(x)$ must be relatively prime. 

\item The degree of $F$ must divide the degrees of each $L_i, i=1..n$.


\end{enumerate}




\section{Summary of the \pFq approach for third order linear ODEs}

The idea consists of assuming that the given linear ODE is one of \pFq
equations \eq{seeds} transformed using \eq{tr} for some $F(x)$ rational in
$x$ and $P(x)$ arbitrary and for some values of the pFq parameters. Resolving
the equivalence is about determining the $F(x)$, $P(x)$ and the values of the
\pFq parameters $\{\alpha, \beta, \gamma, \delta, \eta \}$ such that the
equivalence exists.
An itemized description of the decision procedure to resolve this equivalence, following the presentation the previous sections, is as
follows.

\begin{enumerate}

\item Rewrite the given equation \eq{LinearODE} we want to solve, in normal
form

\begin{equation}
y''' = {\tilde I}_1(x)\, y' + {\tilde I}_0(x)\, y
\ee

\ni where the invariants ${\tilde I}_n(x)$ are constructed using the formulas \eq{invariants}.

\item
\label{HYPER3}
Verify whether an equivalence of the form  $\{x \rightarrow
{(a\,x^k+b)}/{(c\,x^k+d)},\  y \rightarrow P(x)\,y\}$ exists:

\begin{enumerate}

\item Compute ${\tilde J}_n(x)$, the shifted invariants \eq{J_definition}, and use
transformations $x\rightarrow x^k$ to reduce to the integer minimal values
the powers entering the numerator and denominator; i.e., compute $k$ and
$J_n(x)$ in \eq{J_transformation_rule}.

\item Determine the singularities of the $J_n(x)$ and use the classification
of singularities mentioned in section \ref{equiv_section} to tell whether
an equivalence under M\"obius transformations to any of the \3F2, \2F2, \1F2
or \0F2 equations \eq{seeds} exists. 

\item When the equivalence exists, from the singularities of the two
$J_n(x)$ compute the parameters $\{a, b, c, d\}$ entering the M\"obius
transformation \eq{mobius} as well as the hypergeometric parameters $\{\alpha, \beta,
\gamma, \delta, \eta \}$ entering the \PFQ equation \eq{seeds}.

\item Compose the three transformations to obtain one of the form

$$
x \rightarrow \fr{\alpha x^k+\beta}{\gamma x^k+\delta},\ \ y \rightarrow P(x)\,y
$$

mapping the \PFQ equation involved into the ODE being solved.

\end{enumerate}

\item
\label{GCD}
When the equivalence of the previous step does not exist, perform step 1 in the itemization
of section \ref{hyper3}, that is, compute the absolute invariants $L_i$ \eq{Ls} and compute a rational transformation $R(x)$ minimizing
the degrees of the invariants \eq{Ls} of the given
equation

\begin{enumerate}

\item
\label{iteration_step}
When $R(x)$ is not of M\"obius form, change $x \rightarrow R(x)$ rewriting the given equation in canonical form 
and re-enter step \eq{HYPER3} with it, to
resolve the remaining M\"obius transformation and determining the values of the
\pFq parameters.

\end{enumerate}

\item When either of the equivalences considered in steps \eq{HYPER3} or \eq{GCD} exist, compose all the
transformations used and apply the composition to the known solution of the \PFQ
equation to which the equivalence was resolved, obtaining the solution to the
given ODE.

\end{enumerate}

\section{Special cases and MeijerG functions}

Giving a look at the series expansion of any of the \3F2, \2F2, \1F2 or \0F2
functions one can see that there are some different situations that require
special attention at the time of constructing the three independent solutions
to \eq{LinearODE}. Consider for instance the standard \0F2 equation and its three independent solutions,

\be{0F2_and_solution}
\begin{array}{c}
\displaystyle
\yyy + {\frac { \l( \alpha+\beta+1 \r) }{x}}\,\yt
	+ {\frac {\alpha\,\beta}{{x}^{2}}}\,\y1
	- {\frac {1}{{x}^{2}}}\,y = 0
\\*[0.16in]
y={\0F2(\ ;\,\alpha,\beta;\,x)}\,C_1
+ {x}^{1-\beta}{\0F2(\ ;\,2-\beta,1+\alpha-\beta;\,x)}\,C_2
+ {x}^{1-\alpha}{\0F2(\ ;\,2-\alpha,1-\alpha+\beta;\,x)}\,C_3
\end{array}
\ee

\ni where the $C_i$ are arbitrary constants. Expanding in series the first \0F2 function entering this solution we get

\be{0F2_1}
1
+ {\frac {1}{\alpha\,\beta}}x
+ {\frac {1}{2\,\alpha\,\beta\, \l( \alpha+1 \r)  \l( 1+\beta \r) }}{x}^{2}
+ {\frac {1}{6\,\alpha\,\beta\, \l( \alpha+1 \r)  \l( 1+\beta \r)  \l( \alpha+2 \r)  \l( \beta+2 \r) }}{x}^{3}
+ O\l( {x}^{4} \r)
\ee

\ni This series does not exist when $\alpha$ or $\beta$ are zero or
negative integers, and the same happens when the \pFq parameters
entering any of the other two independent solutions is a non-positive integer.
By inspection, however, one of the three \pFq functions entering the solution
in \eq{0F2_and_solution} always exists, because there are no $\alpha$ and
$\beta$ such that the three \0F2 functions simultaneously contain non-positive
integer parameters.

Consider now the second independent solution, ${x}^{1-\beta}{\0F2(\
;\,2-\beta,1+\alpha-\beta;\,x)}$: when $\beta = 1$ it becomes equal to the first
one and so we have only two independent \pFq solutions. In the same way, when
$\alpha = 1$ the first and third solutions entering \eq{0F2_and_solution} are
the same and when $\alpha = \beta$ the second and third solutions are the
same. And when the two conditions hold, that is $\alpha = \beta = 1$, actually
the three solutions are the same. Notwithstanding, in these cases too one of
the three \0F2 solutions always exists.

The same two type of special cases exist for the \1F2, \2F2 and \3F2 function
solutions and the problem at hand consists of having a way to represent the
three independent solutions to \eq{LinearODE} {\it without} introducing
integrals or iterating reductions of order\footnote{Recall that given two
independent solutions, it is always possible to write the third one in terms
of integrals constructed with the two existing solutions, and in the case of a
single solution it is still possible to reduce the order to a second order
linear equation that may or not be solvable.}. For this purpose, we use a set
of 3 MeijerG functions for each of the four \pFq families that can be used to
replace the missing \pFq solutions in these special cases. The key observation
is that at these special values of the last two parameters of the \pFq
functions the MeijerG replacements exist, satisfy the same differential
equation and are independent of the available \pFq function solutions. A table
with these 3 x 4 = 12 MeijerG function replacements is as follows:

\begin{table}[h]
\centering
\label{MeijerGTable}
\caption{MeijerG alternative solutions to the \pFq equations}
\begin{tabular}{|c|c|c|c|}
\hline
\pFq family &  \multicolumn{3}{|c|}{MeijerG functions} \\
\hline
${\0F2(\ ;\,\alpha,\beta;\,x)}$
& $G^{2, 0}_{0, 3}\l(x, \Big\vert\,^{}_{0, 1-\alpha, 1-\beta}\r)$ & $G^{3, 0}_{0, 3}\l(-x, \Big\vert\,^{}_{0, 1-\alpha, 1-\beta}\r)$ & $G^{2, 0}_{0, 3}\l(x, \Big\vert\,^{}_{1-\alpha, 1-\beta, 0}\r)$
\\*[0.1in]
${\1F2(\alpha;\,\beta,\gamma;\,x)}$
& $G^{2, 1}_{1, 3}\l(x, \Big\vert\,^{1-\alpha}_{0, 1-\beta, 1-\gamma}\r)$ & $G^{3, 1}_{1, 3}\l(-x, \Big\vert\,^{1-\alpha}_{0, 1-\gamma, 1-\beta}\r)$ & $G^{2, 1}_{1, 3}\l(x, \Big\vert\,^{1-\alpha}_{1-\gamma, 1-\beta, 0}\r)$
\\*[0.1in]
${\2F2(\alpha,\beta;\,\delta,\gamma;\,x)}$
& $G^{2, 2}_{2, 3}\l(x, \Big\vert\,^{1-\beta, 1-\alpha}_{0, 1-\gamma, 1-\delta}\r)$ & $G^{3, 2}_{2, 3}\l(-x, \Big\vert\,^{1-\beta, 1-\alpha}_{0, 1-\gamma, 1-\delta}\r)$ & $G^{2, 2}_{2, 3}\l(x, \Big\vert\,^{1-\beta, 1-\alpha}_{1-\gamma, 1-\delta, 0}\r)$
\\*[0.1in]
${\3F2(\alpha,\beta,\gamma;\,\delta,\eta;\,x)}$
& $G^{2, 3}_{3, 3}\l(x, \Big\vert\,^{1-\beta, 1-\alpha, 1-\gamma}_{0, 1-\delta, 1-\eta}\r)$ & $G^{3, 3}_{3, 3}\l(-x, \Big\vert\,^{1-\beta, 1-\alpha, 1-\gamma}_{0, 1-\delta, 1-\eta}\r)$ & $G^{2, 3}_{3, 3}\l(x, \Big\vert\,^{1-\beta, 1-\alpha, 1-\gamma}_{1-\delta, 1-\eta, 0}\r)$
\\
\hline
\end{tabular}
\end{table}

\section{Examples}
\medskip
\ni{\bf Equivalence under power composed with M\"obius transformations for the \0F2 class}
\medskip

Consider the third order linear ODE

\begin{eqnarray}
\label{e1}
\lefteqn{\yyy  =
{\frac { \l( 37+2\,\mu+6\,\nu -108\,{x}^{2}\r) }{12\,x \l( x+1 \r)  \l( x-1 \r) }}\,\yt
}
& & 
\\*[0.15in]
\nonumber
& &
+ {\frac { \l( 2\, \l( \nu+6 \r)  \l( 11/2 - \mu\r) + \l( 36\,\nu+294+12\,\mu \r) {x}^{2} -360\,{x}^{4} \r) }{ 24\,{x}^{2} \l( x+1 \r) ^{2} \l( x-1 \r) ^{2}}}\,\y1
- {\frac {16}{x \l( x+1 \r) ^{4} \l( x-1 \r) ^{4}}}\,y
\end{eqnarray}

\ni This equation has two regular singularities at $\{0, \infty\}$ and two
irregular singularities at $\{-1, 1\}$. Following the steps mentioned in the
Summary, we rewrite the equation in normal form and, in step 2.(a), compute the
value of $k$ leading to an equation with minimal degrees entering $J_n(x)$ in
\eq{J_transformation_rule}. The value of $k$ found is $k = 2$ so the equation
from which \eq{e1} is derived changing $x \rightarrow x^2$ is

\begin{eqnarray}
\label{e11}
\lefteqn{
\yyy = {\frac { \l( 6\,\nu+2\,\mu+73 -144\,x\r)}{24\,x \l( x-1 \r) }}\, \yt
 }
& &
\\*[0.1in]
\nonumber
& &
- {\frac { \l( 2\, \l( \nu+8 \r)  \l( \mu + 1/2 \r) - \l( 48\,\nu+16\,\mu+584 \r) x + 576\,{x}^{2} \r) }{ 96\,{x}^{2}\l( x-1 \r) ^{2}}}\,\y1
-{\frac {2}{{x}^{2} \l( x-1 \r) ^{4}}}\,y
\end{eqnarray}

\ni and has invariants with minimal degrees with respect to power
transformations. In step 2.(b), analyzing the structure of singularities of
\eq{e11} we find one regular singularity at the origin and one irregular at
$\infty$.
Using the classification discussed in section 3.2 based on the
degrees with respect to $x$ of the numerators and denominators of the invariants of \eq{e11} as well as
the factors entering these denominators the equation
is identified as equivalent to the \0F2 class under M\"obius transformations
\eq{mobius}. So we proceed with step 2.(c), constructing the M\"obius
transformation and computing the values of the hypergeometric parameters
$\{\mu, \nu\}$ entering the \0F2 equation in \eq{seeds} such that the equivalence
under M\"obius exists, obtaining:

\begin{equation}
\alpha = \nu/4 + 2,\ \ \ \beta = \mu/12 + 1/24,\ \ \ M := x \rightarrow \fr{2\,x}{x-1}
\ee

\ni Composing $M$ above with the power transformation used to obtain \eq{e11} and using the
values above for $\alpha$ and $\beta$, in step 4 we obtain the solution of \eq{e1}

\begin{eqnarray}
\nonumber
y(x) & =  &
{\0F2(\ ;\,\nu/4+2, \mu/12+1/24;\,2\,{\frac {{x}^{2}}{{x}^{2}-1}})}\, C_1
\\*[0.1in]
& & 
+{x}^{-(2+\nu/2)} \l( {x}^{2}-1 \r)^{(1+\nu/4)}{\0F2(\ ;\,-\nu/4, \mu/12-\nu/4-23/24;\,{\frac {2\,{x}^{2}}{{x}^{2}-1}})}\,C_2
\\*[0.1in]
\nonumber
& & + {x}^{({{23}/{12}}-\mu/6)} \l( {x}^{2}-1 \r)^{(\mu/12-{ {23}/{24}})}{\0F2(\ ;\,{{47}/{24}}-\mu/12,\,{{71}/{24}}-\mu/12 + \nu/4;\,{\frac {2\,{x}^{2}}{{x}^{2}-1}})}\,C_3
\end{eqnarray}

\medskip
\ni{\bf Meijerg functions and equivalence under rational transformations for the \1F2 class}
\medskip

Consider the following equation, with no symbolic parameters and only integer
powers

\begin{eqnarray}
\label{k215}
\yyy & = & - {\frac { \l( 6+12\,x -15\,{x}^{2}-6\,{x}^{3}\r) }{x \l( 1+x-{x}^{2} \r)  \l( x+2 \r) }}\,\yt
\\*[0.1in]
& & 
\nonumber
+ {\frac { \l(16 + 48\,x + 36\,{x}^{2} - 20\,{x}^{3} + 9\,{x}^{4} + 81\,{x}^{5} -20\,{x}^{6} - 30\,{x}^{7} - 6\,{x}^{8} \r) }{{x}^{4} \l( x+2 \r) ^{2} \l( 1+x-{x}^{2} \r) ^{2}}}\, \y1
\\*[0.03in]
\nonumber
& & 
- {\frac { \l( x+2 \r) ^{3}}{ \l( 1+x-{x}^{2} \r) ^{2}{x}^{5}}}\,y
\end{eqnarray}

\ni Following steps 1 and 2 in the Summary, we confirm that there exists
no equivalence under \eq{tr_hyper3}, so in step 3 we search
for a rational transformation minimizing the degrees of the invariants \eq{Ls},
finding

\begin{equation}
R(x) = x^2/(1+x)
\ee

\ni Therefore \eq{k215} can be obtained by changing variables $x \rightarrow
R(x)$ in

\be{Jk215}
\yyy =
- {\frac { \l( 3 - 9\,x + 6\,{x}^{2}\r) }{x \l( x-1 \r) ^{2}}}\,\yt
+ {\frac { \l( 1 - 2\,x + 6\,{x}^{2} - 6\,{x}^{3} \r) }{{x}^{3} \l( x-1 \r) ^{2}}}\,\y1
-{\frac {1}{ \l( x-1 \r) ^{2}{x}^{4}}}\,y
\ee

\ni This equation\footnote{The
$c_i$ entering \eq{Jk215} are computed from the minimized $L_j$ by inverting \eq{invariants} and using \eq{Is}.}
thus has invariants with minimal degrees, and has one regular
singularity at 1 and one irregular at the origin. According to the classification in terms of singularities
\eq{Jk215} admits an
equivalence under M\"obius transformations to the \pFq equations (\1F2 case) and hence is solved in the
iteration step 3.(a) mentioned in the summary. When constructing the \pFq solutions to \eq{Jk215}, however,
we find that the
\1F2 parameters in the second list are both equal to 1, so only one \1F2 solution is available, and hence
two of the MeijerG alternative solutions presented in
the table \eq{MeijerGTable} are necessary, resulting in

\begin{equation}
y={{\mbox{$_0${F}$_1$}}(\ ;\,1;\,{\frac {1+x-{x}^{2}}{{x}^{2}}})}\,C_1
+ G^{2, 0}_{0, 2}\l({\frac {1+x-{x}^{2}}{{x}^{2}}}, \Big\vert\,^{}_{0, 0}\r)\,C_2
+ G^{3, 1}_{1, 3}\l({\frac {{x}^{2}-x-1}{{x}^{2}}}, \Big\vert\,^{0}_{0, 0, 0}\r)\,C_3
\ee

\ni Note that the first \pFq function is a {\mbox{$_0${F}$_1$}}. This is due
to the automatic simplification of order that happens when identical parameters are present in both lists of a \1F2
function; this {\mbox{$_0${F}$_1$}} can also be expressed in terms of
Bessel functions.

\section*{Conclusions}

In this work we presented a decision procedure for third order linear ODEs for
computing three independent solutions even when they are not Liouvillian or
when the hypergeometric parameters involved are such that only two or one \pFq
solution around the origin exists. This algorithm solves
complete ODE families we didn't know how to solve before. 

The strategy used is that of resolving an equivalence problem to the \3F2,
\2F2, \1F2 and \0F2 equations, and in doing so, two important generalizations
of the algorithm presented in \cite{hyper3} were developed. First, 
the classification according to singularities and the use of power
composed with M\"obius transformations, presented in \cite{hyper3} for 2nd
order equations, was generalized for third order ones. Second, the idea of
resolving the equivalence mapping into an equation with
invariants with ``minimal degrees under power transformations" was generalized
by determining a transformation mapping into an equation having invariants with
``minimal degrees under general rational transformations". This permits
resolving a much larger class of \pFq equations, defined by changing variables
in \eq{seeds} using $\{x \rightarrow R(x),\, y \rightarrow P(x)\,y\}$ where $R(x)$
is a rational function. Symbolic computation routines implementing this algorithm were integrated into the Maple system in
2007.


Since at the core of the algorithm being presented there is the concept of singularities, two
natural extensions of this work consist of applying the same ideas to
compute solutions for linear ODEs of arbitrary order, where the equivalence
can be solved exactly \cite{MSW_4ODE_challenges}, and for second order
equations under rational transformations, perhaps generalizing the work by
M.Bronstein \cite{bronstein} with regards to {\mbox{$_1${F}$_1$}} solutions to compute also {\mbox{$_2${F}$_1$}} solutions.
Related work is in progress.

%

\end{document}